\documentclass[a4paper,twoside,10pt]{article}
\usepackage{amsthm,graphicx,publaob}

\def\papertype{Contributed paper}

\begin{document}

\title{INVERSE PROBLEM, PRECESSING ORBITS AND AUTONOMOUS FORCES}

\authors{C. BLAGA$^1$}

\address{$^1$Babe\c{s}-Bolyai University, Kog\u{a}lniceanu Street 1, 400084 Cluj-Napoca, Romania}
\Email{cpblaga}{math.ubbcluj}{ro}


\abstract{In the framework of the inverse problem of Dynamics we investigate the
compatibility between a two-parametric family of precessing orbits and autonomous force
fields. In a previous work, for a specific choice of the two parameters in the precessing
family of orbits, we obtained that if the family is compatible with a central force, than
this is derived from a Manev's type potential. In this work we consider another choice of
the parameters in the two-parametric family of precessing orbits, for which we know that
there is no central force compatible and study its compatibility with an autonomous force
field.}

\section{INTRODUCTION}

According to Kepler's first law, the orbits of planets around the Sun should be ellipses,
having the Sun at one of the foci. This law does not take into account the influences of
the other bodies, producing a slow rotation  of these ellipses in their planes. Even so,
there are small differences between the observed precession of planetary orbits and the
one computed based on the Newtonian mechanics. For Mercury, the difference is of about
$43''$/century.

Using the \emph{Inverse Problem Theory} we search the potential responsible for producing
such orbits. In this framework, Xanthopoulos and Bozis (see Xanthopoulos and Bozis 1983)
showed that the potential that produces the family of co-planar ellipses with a common
focus and arbitrary eccentricity, magnitude and orientation of their major axes has to be
the Newtonian potential, in other words, they were able to show that it is possible to
deduce Newton's gravitational law using only the first Kepler's law (applied, to several
planets, not only to a single one).

The equation of a precessing conic is
\begin{equation}\label{1}
r(\theta)=\frac{p}{1+e\cos b(\theta-\theta_0)}
\end{equation}
where $p=a(1-e^2)$ is the the semilatus rectum, $a$ the semimajor axis, $e$ the
eccentricity of the conic section and $\theta_0$ the orientation of the semimajor axis.
The rotation of the semimajor axis in the plane of the the conic section is given by $b$.

We write this equation in the form
\begin{equation}\label{2}
f(r,\theta,\theta_0)\equiv \frac{p-r}{r\cos b(\theta-\theta_0)}=e
\end{equation}
and consider it as the equation of a two-parametric family of curves, the two parameters
being the orientation of the semimajor axis, $\theta_0$ and the eccentricity $e$, while
the other two parameters ($p$ and $b$) are assumed to be functions of the first ones.

\section{PRECESSING ORBITS AND CENTRAL FORCES}

The compatibility between a two-parametric family of orbits and a central force was
discussed  by Borghero, Bozis and Melis (see Borghero, Bozis and Melis 1999), paper in
which they obtained the conditions verified by a two-parametric family of orbits whose
members are trajectories for  masspoints moving in a central field. These criteria are
fulfilled for the family~(\ref{2}) and the force compatible to it is even a conservative
one (see Blaga 2005). Its potential turns out to be the potential
\begin{equation}\label{3}
V(r)=-F_0\left[\frac{b^2}{r}+\frac{p(1-b)^2}{2r^2}\right],
\end{equation}
known as Manev's potential.

\theoremstyle{remark}
\newtheorem*{rem}{Remark}
\begin{rem}
For $b=1$ the family of orbits~(\ref{1}) becomes a family of orbits without precession
and the potential~(\ref{3}) reduces to the Newtonian one.

If we choose $b$ or $p$ as parameter in~(\ref{2}) we get a negative answer to
the question concerning the compatibility of this two parametric family of
orbits with a central force. There is \emph{no} central force that admits as
orbits all the members of the two parametric family
\begin{equation}\label{4}
f(r,\theta,b)\equiv r [1+ e \cos b(\theta-\theta_0) ]= p
\end{equation}
with $b$ and $p$ parameters, and fixed $e$ and $\theta_0$.
\end{rem}

\section{TWO-PARAMETRIC FAMILIES OF ORBITS AND AUTONOMOUS FORCES}

But are there any autonomous dynamical systems which generates this two-parametric family
of orbits? The compatibility between a two-parametric family of orbits and an autonomous
conservative dynamical system was discussed by Bozis (Bozis 1983), paper in which the
author gave criteria to test the existence of the solution and, in the case of a positive
answer, a method to compute it.

We consider the family of planar precessing orbits in Cartesian coordinates
\begin{equation}\label{5}
f(x,y,b)\equiv (x^2+y^2) [1+ e \cos b(\theta-\theta_0) ]^2= p^2
\end{equation}
with $\theta = \arctan(y/x)$, $b$ and $p$ parameters, and fixed $e$ and $\theta_0$.

If there is a dynamical system which generates this family of orbits, then the
force components $X(x,y)$ and $Y(x,y)$ are related to the family of orbits
through:
\begin{equation}\label{6}
-X_x+\frac{1}{\gamma} X_y- \gamma Y_x +Y_y = \lambda X + \mu Y \,
\end{equation}
where the coefficients of the equation are
\begin{equation}\label{7}
\lambda = \left( -\Gamma_x + \frac{1}{\gamma} \Gamma_y \right) \Gamma^{-1} \qquad
\mbox{and} \qquad \mu = \lambda \gamma + \frac{3 \Gamma}{\gamma},
\end{equation}
and
\begin{equation}\label{8}
\gamma = \frac{f_y}{f_x} \qquad \mbox{and} \qquad \Gamma = \gamma \gamma_x - \gamma_y.
\end{equation}

If the dynamical system is a conservative one, then $X_y=Y_x$.

The coefficients of the equation~(\ref{6}) are, generally, functions of $x$,
$y$ and $b$, but we seek solutions $X(x,y)$ and $Y(x,y)$ independent of $b$, in
other words
\begin{equation}\label{9}
X_b=Y_b=0 \, .
\end{equation}

Using this condition, Bozis (see Bozis 1983) classified the cases that can
arise and gave an algorithm to find the solution if it exists. We need to
introduce the functions
\begin{equation}\label{10}
L=-\frac{\gamma^2}{(1+\gamma^2) \gamma_b} \lambda_b, \qquad \mbox{and} \qquad
M=-\frac{\gamma^2}{(1+\gamma^2) \gamma_b} \mu_b \, .
\end{equation}

The compatibility of the family of orbits~(\ref{5}) with a conservative autonomous
dynamical system is ruled by the following

\theoremstyle{plain}
\newtheorem*{prop}{Proposition}
\begin{prop}
If
\begin{equation}\label{11}
L_b \neq 0 \, , \qquad M_b \neq 0 \, , \qquad \left( \frac{L}{M} \right)_b \neq 0,
\end{equation}
and the conditions $\left( \frac{L_b}{M_b} \right)_b = 0$ \mbox{and} $\left(
\frac{L+M \rho - \rho_x}{\rho} \right)_y = (L+M \rho)_x$, with $L$ and $M$
given by~(\ref{10}), and $\rho=-\frac{L_b}{M_b}$, are fulfilled, then the
problem admits a solution which is found solving the equations:
\begin{equation}\label{12}
X_x=\left( -\frac{D}{L_b} - \frac{\rho_x}{\rho} \right) X \, , \quad X_y=\frac{D}{M_b} X
\, , \quad Y= \rho X\, ,
\end{equation}
where $D=L M_b-M L_b$. If one of the above conditions is not satisfied, then no solution
exists.
\end{prop}

\section{THE TWO-PARAMETRIC FAMILY OF PRECESSING ORBITS}

For the family~(\ref{4}) we find
\begin{equation}\label{13}
\gamma = \tan (\theta-\omega)
\end{equation}
and
\begin{equation}\label{14}
\Gamma = - \frac{\cos^2 \omega}{\sqrt{x^2+y^2} \, \cos^2 (\theta -\omega)} \, \cdot
\,\frac{1+ e \cos \theta'- b^2 e \cos \theta'}{cos \theta + e \cos \theta \cos \theta' +
b e \sin \theta \sin \theta'}
\end{equation}
where
\begin{equation}\label{15}
\tan \omega = \frac{b e \sin \theta'}{1+ e \cos \theta'} \quad \mbox{and} \quad \theta'=
b (\theta - \theta_0) \,.
\end{equation}

To simplify the computation of the functions needed in our analysis, we observe that the
function $f$ from ~(\ref{4}) is  homogenous of degree one in $x$ and $y$. The function
$\gamma = f_y/f_x$ is homogeneous of degree zero in $x$ and $y$ and we can consider it as
a function in one variable $z=y/x$. Denoting by $\dot{\gamma}$ and $\ddot{\gamma}$ the
first and second derivative with respect to $z$, we obtain that
\begin{equation}\label{16}
\lambda = \frac{(\gamma z + 1) \ddot{\gamma} + z \, \dot{\gamma} ^ 2 + 2 \gamma \,
\dot{\gamma}}{x \, \gamma \, \dot{\gamma}} \quad \mbox{and} \quad \mu = \frac{(\gamma z +
1) (\ddot{\gamma} \gamma - 3 \dot{\gamma} ^2) + z \, \gamma \, \dot{\gamma} ^ 2 + 2
\gamma^2 \, \dot{\gamma}}{x \, \gamma \, \dot{\gamma}}
\end{equation}

\section{RESULTS AND CONCLUSIONS}

After tedious but straightforward calculations we find that
\begin{equation}\label{17}
\frac{d \gamma}{d z} = \frac{\cos^2 \theta \, \cos \omega}{\cos (\theta - \omega)} \,
\cdot \, \frac{1 + e \, (1-b^2) \, \cos \theta'}{( \cos \theta + e \, \cos \theta \, \cos
\theta' + e \, b\, \sin \theta \, \sin \theta')}
\end{equation}
and
\begin{equation}\label{18}
\frac{d^2 \gamma}{d z^2} = \frac{e \, b \, \cos^3 \theta}{\cos ( \theta - \omega)} \,
\cdot \, \frac{c_3 b^3 + c_2 b^2 +c_1 b + c_0}{( \cos \theta + e \, \cos \theta \, \cos
\theta' + e \, b\, \sin \theta \, \sin \theta')^2 }
\end{equation}
 where
\begin{equation}
c_3 = e \, \cos \omega \, \sin 2 \theta \, , \qquad c_2 =e \sin 2\theta \, \cos \omega +
\sin \theta' \, \cos \theta \, \cos (\theta + \omega) \, ,
\end{equation}
\begin{equation}
c_1 =- \sin 2 \theta \,\cos \omega \, \left( \cos \theta' +e \right) \, , \quad c_0 = -2
\, \cos \omega \sin \theta'  \left( e\, \cos \theta'
 +1 \right) \, .
\end{equation}

Replacing these functions in~(\ref{16}) we find that the functions $\lambda$, $\mu$, $L$
and $M$ depend on $b$ directly and through $\theta'$ and $\omega$. It is now
straightforward to find that $L_b$, $M_b$, $(L/M)_b$ and $(L_b/M_b)_b$, functions which
are not zero. According to the proposition of section 3 these results imply that there is
\emph{no} autonomous conservative dynamical system compatible with the two-parametric
family of orbits~(\ref{4}).

\begin{center}
\bf{ Acknowledgement}
\end{center}
During the preparation of my talk for this conference I  had discussions with dr. George
Bozis, from the University of Thessaloniki. I would like to express my thanks for all the
help he gave me.

\references

Blaga C. : 2005, \journal{Prof. G. Manev's Legacy in Contemporary Astronomy, Theoretical
and Gravitational Physics, eds. Gerdjikov V., Tsvetkov M.}, Sofia: Heron Press 134.

Borghero F., Bozis G., Melis A. : 1999, \journal{Rendiconti di Matematica}, \vol{19},
303.

Bozis G. : 1995, \journal{Inverse Problems Journ.}, \vol{11}, 687.

Bozis G. : 1983, \journal{Celestial Mechanics}, \vol{31}, 129.

Xanthopoulos B. and Bozis G. : 1983, \journal{Dynamical Trapping and Evolution in the
Solar System, eds. Y Kozai and V Markellos}, IAU Coll. 74, Dordrecht: Reidel, 353.

\endreferences

\end{document}